\documentclass{amsart}
\usepackage{amsfonts,amssymb,amsmath,amsthm}

\usepackage{amsmath,amssymb,amsfonts}%
\usepackage{amsthm}

\makeatletter
\@namedef{subjclassname@2020}{\textup{2020} Mathematics Subject Classification}
\makeatother

\newtheorem{thrm}{Theorem}
\newtheorem{prop}[thrm]{Proposition}%
\newtheorem{ex}{Example}%
\newtheorem{remark}{Remark}%
\newtheorem{cor}[thrm]{Corollary}
\newtheorem{definition}{Definition}%
\newtheorem{lem}[thrm]{Lemma}
\newtheorem{Question}[thrm]{Problem}

\begin{document}
\title{On the universal approximation of real functions with varying domain}


\author[W. Jung]{W. Jung}
\address{Healthcare Datascience Center, Konyang University Hospital,
Daejeon, South Korea.}
\email{wcjungdynam@gmail.com}
 
\author[C.A. Morales]{C.A. Morales}
\address{Beijing Advanced Innovation Center for Future Blockchain and Privacy Computing, Beihang University, Beijing, China \& Beijing Academy of Blockchain and edge Computing, Beijing, 100086.}
\email{morales@impa.br}

\author[L.T.T. Tran]{L.T.T. Tran}
\address{BNF Technology Techno2-ro 170-10, Yuseong District, Daejeon, South Korea.}
\email{tolinh.tran91@gmail.com}

\begin{abstract}
We establish sufficient conditions for the density of shallow neural networks \cite{C89} on the family of continuous real functions defined on a compact metric space, taking into account variations in the function domains. For this we use the Gromov-Hausdorff distance defined in \cite{5G}.
\end{abstract}

\keywords{Approximation, Metric space, Gromov-Hausdorff distance}
\subjclass[2020]{Primary 41A65; Secondary 54C2}

\maketitle




\section{Introduction}

\noindent
The goal of this paper is to find sufficient conditions for the shallow neural network generated by a
family of real-valued continuous functions on compact metric spaces to satisfy a “universal approximation property”. We will not assume that the domains of the functions coincide.
To make this precise, we use the Gromov-Hausdorff distance on the space of compact metric spaces \cite{bbi} and its related distance for functions \cite{5G}.
Our main result (Theorem \ref{pertinho}) provides the desired sufficient conditions. A concrete corollary (Corollary \ref{pertitinho}) asserts that the shallow neural network generated by functions defined on finite metric spaces is dense in the space of all continuous functions on compact metric spaces. 
This is new and stronger than the well-known fact that the finite metric spaces are dense in the space of compact metric spaces endowed with the Gromov-Hausdorff metric \cite{bbi} (see for instance \cite{ju} or the proof of Example \ref{bomber}).

A potential application of our result is in the representation of shapes by point clouds.
This is because finite spaces can be used to approximate objects with the
Gromov-Hausdorff distance measuring the dissimilarity of these objects
(see for instance Section 1 in \cite{5G}). 
Then, real-valued functions over finite spaces can therefore represent some shape data to be measured.
Let us state our results in a precise way.

Let $X$ be a metric space.
Denote
$\mathcal{C}(X)$ the set of continuous maps $f:X\to \mathbb{R}$. This is a
Banach algebra under the usual operations and the supremum norm
$$
\|f\|_\infty=\sup_{x\in X}|f(x)|.
$$
Define the {\em $C^0$-distance} $d_{C^0}(f,g)=\|f-g\|_\infty
$, for all $f,g\in \mathcal{C}(X).$
An {\em activation function} is a
continuous function $\sigma:\mathbb{R}\to \mathbb{R}$.
It is of {\em sigmoidal type} if
$$
\lim_{t\to-\infty}\sigma(t)=0\quad\mbox{ and }\quad
\lim_{t\to\infty}\sigma(t)=1.
$$
Following the standard terminology of neural networks with a single hidden layer we introduce the definition below.

\medskip
 
\begin{definition}
The {\em shallow neural network} $N_\sigma(F)$ generated by an activation function $\sigma$ and $F\subset \mathcal{C}(X)$
is $\emptyset$ (if $F=\emptyset$) or the collection of functions $g:X\to \mathbb{R}$ of the form
$$
 g(x)=\sum_{j=1}^Na_j\sigma(f_j(x)+\theta_j)\quad\quad (\forall x\in X)
$$
where $N$, $a_1,\cdots, a_N,\theta_1,\cdots, \theta_N$ and
$f_1,\cdots, f_N$ run over $\mathbb{N}$, $\mathbb{R}$, and $F$ respectively.
 \end{definition}

\medskip

Cybenko \cite{C89} proved that, for any $n\in\mathbb{N}$, 
the shallow neural network generated by an activation function of sigmoidal type and the dual of $\mathbb{R}^n$ on the $n$-dimensional cube $[0,1]^n$ is dense in $C([0,1]^n)$ w.r.t the $d_{C^0}$-distance. 
Our goal is to generalize this result by allowing the domain of functions to vary in the Gromov-Hausdorff space.

To make this precise, let $\mathcal{M}$ be the set of all compact metric spaces.
Define
$$
\mathcal{C}=\bigcup_{X\in \mathcal{M}}\mathcal{C}(X).
$$
Let $D:\mathcal{C}\to \mathcal{M}$ be the {\em domain map} defined by
$$
D(f)=X\quad\iff\quad X\mbox{ is the domain of }f.
$$
Given $X,Y\in \mathcal{M}$ and $\epsilon>0$ we say that
$i:X\to Y$ is an {\em $\epsilon$-isometry} if
\begin{itemize}
\item
$|d(i(x),i(x')-d(x,x')|\leq\epsilon$ ($\forall x,x'\in X$);
\item
For every $y\in Y$ there is $x\in X$ such that
$d(i(x),y)\leq\epsilon$.
\end{itemize}

The {\em Gromov-Hausdorff distance} in $\mathcal{M}$ is defined by
$$
d_{GH}(X,Y)=\inf\{\epsilon>0:\exists\mbox{ $\epsilon$-isometry from $X$ to $Y$ and vice versa}\}.
$$
The {\em $C^0$-Gromov-Hausdorff distance} in $\mathcal{C}$ is defined by
\begin{eqnarray*}
d_{GH^0}(f,g)&=&\inf\{\epsilon>0:\exists\mbox{ $\epsilon$-isometries $i:D(f)\to D(g)$}\\
& & \mbox{ and $j:D(g)\to D(f)$ s.t. }
\|g\circ i-f\|_\infty\leq\epsilon\mbox{ and } \\
& & \quad\quad\quad \|g-f\circ j\|_\infty\leq\epsilon\}.
\end{eqnarray*}
(The latter is not the original definition \cite{5G} but an equivalent one based on the approach in \cite{AM17}).
Both distances generate topologies induced by the open balls.
For brevity, we say that a
given collection of subsets of $\mathcal{C}$ is {\em dense in $\mathcal{C}$}
if it is dense in $\mathcal{C}$ with respect to $d_{GH^0}$.

For every $\mathcal{F}\subset \mathcal{C}$ we define
$\mathcal{F}:\mathcal{M}\to\mathcal{C}$ by
$$
\mathcal{F}(X)=\mathcal{F}\cap \mathcal{C}(X)\quad\quad\forall X\in \mathcal{M}.
$$

\begin{definition}
The {\em shallow neural network}  generated by an activation function $\sigma$ and $\mathcal{F}\subset \mathcal{C}$
is defined by
$$
\mathcal{N}_\sigma(\mathcal{F})=\bigcup_{X\in\mathcal{M}}N_\sigma(\mathcal{F}(X)).
$$
\end{definition}

The main problem is the following one.

\medskip

\begin{Question}
Find sufficient conditions for the shallow neural network generated by an activation function $\sigma$ and a collection $\mathcal{F}\subset \mathcal{C}$ to be dense in $\mathcal{C}$.
\end{Question}

\medskip

To attack this problem we first we note that both $D(\mathcal{N}_\sigma(\mathcal{F}))=D(\mathcal{F})$
and $\mathcal{F}\subset D^{-1}(D(\mathcal{F}))$.
Then, for $\mathcal{N}_\sigma(\mathcal{F})$ to be dense in $\mathcal{C}$ it is necessary that
\begin{equation}
\label{rosa}
D^{-1}(D(\mathcal{F})) \text{ is dense in } \mathcal{C}
\end{equation}

\medskip

Clearly, this condition is insufficient (see the next section).
Second, let $B(X)$ be a set of regular signed Borel measures of $X\in \mathcal{M}$.
Given a measurable map $\varphi:X\to Y$ where $Y\in \mathcal{M}$, we define the pushforward $\varphi:B(X)\to B(Y)$ by $\varphi_*(\mu)=\mu\circ \varphi^{-1}$ for $\mu\in B(X)$.
The following definition can be found in p. 112 of \cite{ST09}.

\medskip

\begin{definition}
Given $X\in \mathcal{M}$, we say that $F \subseteq \mathcal{C}(X)$ separates $B(X)$ if
$$
\mu\in B(X)\quad\mbox{ and }\quad
f_*(\mu) = 0 \:\: (\forall f \in F)
\quad\Longrightarrow\quad \mu=0.
$$
\end{definition}

\medskip

Next, we introduce some basic terminology about fibrewise sets (see \cite{j} for details). A {\em fibrewise space} over a topological space $M$ consists of a topological space $C$ together with a continuous function $D:C\to M$, called the projection. For each $x\in M$ the fibre over $x$ is the subset $M_x=D^{-1}(x)$ of $C$;
fibres may be empty since $D$ is not required to be surjective.
When the fibres have a similar structure, the resulting space is termed a fibre bundle, emphasizing the compatibility of the fibres and the continuity of their transitions over \(M\).
We will use the following definition.

\medskip

\begin{definition}
\label{didi}
A fibrewise set $C$ over $M$ is {\em locally sliceable} if for every $x\in M$, $f\in M_x$ and every neighborhood $U$ of $f$ in $C$ there is a neighborhood $W$ of $x$ in $M$ such that $M_{x'}\cap U\neq\emptyset$ for every $x'\in W$.
\end{definition}

\medskip

This definition is weaker but closely related to the one with the same name in Definition 1.16 of \cite{j}.
The fibrewise sets we are going to deal with are those with total space within $\mathcal{C}$, base space within $\mathcal{M}$ and projection giving by the restriction of $D$.
About them we have the following definition.

\medskip

\begin{definition}
We say that $\mathcal{F}\subset \mathcal{C}$ is {\em locally sliceable}
if the fibrewise subset of $\mathcal{C}$ with total space $D^{-1}(D(\mathcal{F}))$ and base $D(\mathcal{F})$ is locally sliceable.
\end{definition}

\medskip

We shall present examples  related to the above definitions in the next section.
Now, we can state our main result.

\medskip
\begin{thrm}
\label{pertinho}
Let $\mathcal{F}$ a locally sliceable subset of $\mathcal{C}$ satisfying \eqref{rosa} and $\sigma$ be an activation function of sigmoidal type. If
$$
\{X\in D(\mathcal{F}):\mathcal{F}(X)\mbox{ separates }B(X)\mbox{ and }\mathbb{R}\mathcal{F}(X)\subset \mathcal{F}(X)\}
$$
is dense in $D(\mathcal{F})$,
then the shallow neural network generated by $\sigma$ and
$\mathcal{F}$ is dense in $\mathcal{C}$.
\end{thrm}

\medskip

An anonymous referee asked if the locally sliceable condition can be removed from this theorem.
However, our methods depend heavily on this condition so we do not know if it can be removed or not.

A short application of this result is as follows.
Let $\mathcal{F}_{fin}$ be the family of all real-valued functions defined on finite metric spaces.
In Example \ref{bomber} we will show that $\mathcal{F}_{fin}$ satisfies the conditions this theorem. Then, the following corollary holds.

\medskip
\begin{cor}
\label{pertitinho}
The shallow neural network generated by an activation function of sigmoidal type
and $\mathcal{F}_{fin}$ is dense in $\mathcal{C}$.
\end{cor} 

\medskip

The remainder of the paper is divided as follows.
In Section \ref{sec2}, we give some examples, properties of $\mathcal{N}_\sigma(\mathcal{F})$ and provide families separating signed Borel measures.
In Section \ref{sec3}, we prove Theorem \ref{pertinho}
by extending Cybenko methods \cite{C89} to the varying domain case.

The authors would like to thank the Vietnam Institute for Advanced Study in Mathematics (VIASM) for their hospitality and financial support.

\section{Examples}
\label{sec2}

\noindent
In this section, we present some examples to support the main definitions of this work.
Unless otherwise stated $\sigma$ denotes an activation function, $X\in \mathcal{M}$ and $F\subset \mathcal{C}(X)$.

\medskip
\begin{ex}
If $\sigma=0$ then $N_\sigma(F)=\{0\}$.
\end{ex}
\medskip

By this example we always assume $\sigma\neq0$.
Given $n\in\mathbb{N}$, we denote by $I_n=[0,1]\times\overset{(n)}{\cdots}\times [0,1]$.
The set of linear operators from $\mathbb{R}^n$ to $\mathbb{R}$ is the dual of $\mathbb{R}^n$ denoted by $(\mathbb{R}^n)^*$.

\medskip
\begin{ex}
\label{bucle}
If $X = I_n$ and $F= \{ L|_X : L \in (\mathbb{R}^n)^* \}$ be the dual of $\mathbb{R}^n$ restricted to $X$, then $N_\sigma(F)$ is the family of functions considered in Eq. (1) p. 303 of \cite{C89}.
\end{ex}
\medskip

Denote by $Con\subset \mathcal{C}$ the set of real-valued constant functions defined on compact metric spaces.

\medskip
\begin{ex}\label{ex2}
Since $\sigma\neq0$, we have $\sigma\circ Con = Con$.
\end{ex}

\medskip

The following example shows that the condition \eqref{rosa} is not sufficient to guarantee the density of $\sigma(\mathcal{F})$ on $\mathcal{C}$.

\medskip

\begin{ex}
The collection of all zero maps $\mathcal{F}=\{0:X\to \mathbb{R}\}_{X\in \mathcal{M}}$ satisfies \eqref{rosa}, but $N_\sigma(\mathcal{F})$
is not necessarily dense in $\mathcal{C}$ for it consists of constant maps
(by Example \ref{ex2}).
\end{ex}

\medskip

Given $F\subset \mathcal{C}(X)$ for some $X\in\mathcal{M}$, we observe that $N_\sigma(F)$ is a subspace but not necessarily a subalgebra of $C(X)$. Therefore, the Stone-Weierstrass theorem may not apply to determine when $N_\sigma(F)$ is dense in $C(X)$. For instance, we can consider the conditions under which $Con(X) \subseteq N_\sigma(F)$. Below, we provide a very satisfactory answer.

\medskip
\begin{remark}
If $F \cap Con(X) \neq \emptyset$, then $Con(X) \subseteq N_\sigma(F).$
\end{remark}
\medskip

Next, we consider the problem when
$\sigma(F)$ separates points. 

\medskip
\begin{remark}
If $\sigma$ is injective and $F$ separate points of $X$, then $N_\sigma(F)$ separates points as well.
\end{remark}
\medskip

\begin{proof}
Pick $x \neq y \in X$ and take $f \in F$ such that $f(x) \neq f(y)$. Take $N=1, a_1=1$ and $\theta=0$, thus $g = \sigma\circ f \in \sigma(F)$ and $g(x) \neq g(y)$.
\end{proof}

Next we consider the problem when the $
N_\sigma(F)$ is an algebra of $\mathcal{C}(X)$.

\medskip
\begin{ex}
Suppose there is a real number $A \neq 0$ such that
\begin{equation}
\label{multi}
\sigma(t) \cdot \sigma(s) = A \sigma(t \cdot s) \:\: {\text for \: all} \:\: t,s \in \mathbb{R}.
\end{equation}
If $F$ is an algebra of $\mathcal{C}(X)$, then so is $N_\sigma(F)$.
\end{ex}
\medskip

\begin{proof}
Given $N,M \in \mathbb{N}$, for every $a_i, b_j, \theta_i, \check{\theta}_j \in \mathbb{R}$, and $f_i, \check{f_i} \in F$,
\begin{eqnarray*}
& & \Bigg[ \sum_{i=1}^{N} a_i \sigma(f_i(x)+\theta_i) \Bigg] \cdot \Bigg[ \sum_{j=1}^{M} b_j \sigma(\check{f}_j(x)+\check{\theta}_j) \Bigg] \\
& = & \sum_{i=1}^{N} \sum_{j=1}^{M} a_i b_j [ \sigma(f_i(x)+\theta_i) \cdot \sigma(\check{f}_j(x)+\check{\theta}_j)] \\
& = & \sum_{i=1}^{N} \sum_{j=1}^{M} a_i b_j A\sigma[(f_i(x)+\theta_i) \cdot (\check{f}_j(x)+\check{\theta}_j)] \\
& = & \sum_{i=1}^{N} \sum_{j=1}^{M} a_i b_j A\sigma[(f_i \cdot \check{f}_j)(x)+ \check{\theta}_j f_{i}(x)+ \theta_i \check{f}_j(x)+ \theta_i \check{\theta}_j].
\end{eqnarray*}
Since $F$ is an algebra, we have
$(f_i \cdot \check{f}_j)+ \check{\theta}_j f_{i}+ \theta_i \check{f}_j\in F$ (for all $i=1,\cdots, N$ and $j=1,\cdots M$)
hence $N_\sigma(F)$ is an algebra.
\end{proof}

Examples of activation functions $\sigma$ satisfying
\eqref{multi} are:

\begin{itemize}
\item[(1)] $\sigma(t) = a |t|$,
\item[(2)] $\sigma(t) = a \cdot t^p$ for some $p \in \mathbb{N}$,
\item[(3)] $\sigma(t) = a |t|^p, \: p \geq 0, \: a \neq 0$.
\end{itemize}

\medskip
\begin{ex}
We say that $F \subseteq C(X)$ is a Stone-Weierstrass family if it is an algebra containing a non-zero constant function and separating points
(hence they are dense in $C(X)$).
If $\sigma$ is injective and satisfies \eqref{multi} and $F$ is a Stone-Weierstrass family, then $N_\sigma(F)$ is also a Stone-Weierstrass family.
\end{ex}
\medskip

Let us give some examples of families separating the signed Borel measures.

\medskip

\begin{ex}
Let $X = I_n$ and $F= \{ L|_X : L \in (\mathbb{R}^n)^* \}$ be the dual of $\mathbb{R}^n$ restricted to $X$. Then,  $F$ separates $B(X)$.
\end{ex}
\medskip

\begin{proof}
(See \cite{C89}).
Let $\mu \in B(X)$ such that $f_*(\mu) = 0$ for all $f \in F$. Then, for every bounded measurable $h : \mathbb{R}^n \to \mathbb{R}^n$ we have
$$
\int_{I_n} h(\langle y,x\rangle) d\mu = 0, \: \text{for all} \: x,y \in \mathbb{R}^n.
$$
In particular, for $h(u) = sin(u)$ and $h(u) = cos(u)$ one has
$$
\int_{I_n} e ^{\sqrt{-1}\langle m,x\rangle} d\mu=\int_{I_n} cos(\langle m,x\rangle)d\mu + \sqrt{-1} \int_{I_n} sin(\langle m,x\rangle) d\mu =0$$
for all $m\in \mathbb{R}^n$ with integer entries. Then, the Fourier transform of $\mu$ is $0$ and so $\mu = 0$ as well. Therefore, $F$ separates $B(X)$.
\end{proof}

\medskip
\begin{ex}
Let $X$ be a compact subset of $\mathbb{R}$.
If $f:X\to \mathbb{R}$ and
$f(x)=x$ for every $x$, then $F = \{f  \}$
separates $B(X)$. 
\end{ex}
\medskip

\begin{proof}
Just note that $f_*(\mu) = \mu$ for all $\mu\in B(X)$, then the conclusion follows.
\end{proof}

Now, we give an example of a family satisfying the hypotheses of our main theorem.

\medskip
\begin{ex}
\label{bomber}
The family $\mathcal{F}_{fin}=\{f\in \mathcal{C}:D(f)\mbox{ is finite}\}$ satisfies the hypothesis of Theorem \ref{pertinho}.
\end{ex}
\medskip

\begin{proof}
Since $D(\mathcal{F}_{fin})$ is the set of finite metric spaces
which dense in $\mathcal{M}$
(see Exercice 7.4.9 in \cite{bbi}), $D(\mathcal{F}_{fin})$ is dense in $\mathcal{M}$.
Therefore, \eqref{rosa} holds.

Now take $\hat{f}\in D^{-1}(D(\mathcal{F}_{fin}))$ and $\Delta>0$.
Choose $0<\rho<\Delta$ such that
$$
d(a,b)<\rho\quad(a,b\in D(\hat{f}))\quad\Longrightarrow\quad |\hat{f}(a)-\hat{f}(b)|<\Delta.
$$
Suppose $\hat{X}\in D(\mathcal{F}_{fin})$ satisfies
$$
d_{GH}(\hat{X},D(\hat{f}))<\frac{\rho}4.
$$
Then, there is a $\frac{\rho}4$-isometry $i:\hat{X}\to D(\hat{f})$.
By Lemma 1.6 in \cite{lm} there is a $\frac{3\rho}{4}$-isometry $j:D(\hat{f})\to\hat{X}$ such that
$$
d_{C^0}(i\circ j,id_{D(\hat{f})})<\frac{\rho}4,
$$
where $id_Z$ denotes the identity of $Z\in \mathcal{M}$.
Define
$$
\check{f}=\hat{f}\circ i.
$$
Then,
$$
\|\check{f}-\hat{f}\circ i\|_\infty=0<\Delta.
$$
Moreover,
$$
|\check{f}\circ j(x)-\hat{f}(x)\|_\infty=|\hat{f}\circ i\circ j(x)-\hat{f}(x)|=|\hat{f}(a)-\hat{f}(b)|<\Delta
$$
for $a=i\circ j(x)$ and $b=x$ satisfies $d(a,b)<\rho$.
Then,
$$
\|\check{f}\circ j-\hat{f}\|_\infty<\Delta.
$$
Since both $i$ and $j$ are $\delta$-isometries we get
$$
d_{GH^0}(\check{f},\hat{f})<\Delta.
$$
Therefore, $D(\mathcal{F}_{fin})$  (and so $\mathcal{F}_{fin}$) are locally sliceable.

Finally, we observe that
if $X\in D(\mathcal{F}_{fin})$, then $X$ is finite so
$\mathcal{F}_{fin}(X)=\mathcal{C}(X)$ thus
$\mathcal{F}_{fin}(X)$ separates $B(X)$ and satisfies
$\mathbb{R}\mathcal{F}_{fin}(X)\subset \mathcal{F}_{fin}(X)$ for all $X\in D(\mathcal{F}_{fin})$.
This proves the result.
\end{proof}

\section{Proof of Theorem \ref{pertinho}}
\label{sec3}

\noindent
First we prove the following result.

\medskip
\begin{prop}\label{MT2}
Let $X$ be a compact metric space, and $\sigma : \mathbb{R} \to \mathbb{R}$ be an activation function of sigmoidal type. If $F \subseteq \mathcal{C}(X)$ separates $B(X)$
and $\mathbb{R}F \subset F$, then $N_\sigma(F)$ is dense in $\mathcal{C}(X)$ with respect to the distance $d_{C^0}$ (or $C^0$-dense for short).
\end{prop}
\medskip

Observe the disparity between this statement and the findings in \cite{B21}, where the assumption entails $\sigma$ being non-polynomial and $F$ separate points to establish the density of $N_\sigma(F)$'s {\em span} in $C(X)$.

The proof this proposition is based on the methods by Cybenko \cite{C89}. We just need some preliminaries. The first one is the definition of discriminatory $\sigma$, motivated by \cite{C89}.

Hereafter $X$ is a compact metric space, $\sigma$ is an activation function, and $F \subseteq \mathcal{C}(X)$.

\medskip
\begin{definition}
We say $\sigma$ is discriminatory with respect to $F$ if
$$\int_X \sigma(f(x) + \theta) d\mu = 0 \:\: \text{for all} \:\: f \in F,\, \theta \in \mathbb{R}
\quad\Longrightarrow\quad \mu=0.$$
\end{definition}
\medskip

The statement of lemmas \ref{lemma1} and \ref{lemma2} are similiarly in \cite{C89} (the proof also) by replacing $\mathbb{R}^n$ by $X$. 

\medskip
\begin{lem}\label{lemma1}
If $\sigma$ is discriminatory with respect to $F$, then $N_\sigma(F)$ is dense in $\mathcal{C}(X)$.
\end{lem}
\medskip

\begin{proof}
We know that $N_\sigma(F)$ is a subspace of $\mathcal{C}(X)$. If $\sigma(F)$ is not dense (i.e $\overline{\sigma(F)} \neq \mathcal{C}(X)$), there is bounded linear functional $L \neq 0$ of $\mathcal{C}(X)$ such that $L|_{\sigma(F)} = 0$.
By Riesz representation theorem, there is $\mu \in B(X)$ such that
$$L(h) = \int_{X} h(x) d \mu(x), \quad\quad \forall h \in \mathcal{C}(X).$$
Then, for every $h \in N_\sigma(F)$ one has
$$\int_{X} h(x) d\mu(x) = 0 .$$
So, $$\int_{X} \sigma(f(x)+\theta) d\mu(x), \quad\quad  \forall f \in \mathcal{F}(X), \theta \in \mathbb{R}$$
thus $\mu =0$ since $\sigma$ is discriminatory with respect to $F$. Therefore, $L=0$ which is absurd. This contradiction completes the proof.
\end{proof}

\medskip
\begin{lem}\label{lemma2}
If $\sigma$ is of sigmoidal type, $F$ separates $B(X)$ and $\mathbb{R}F\subset F$, then  $\sigma$ is discriminatory with respect to $F$.
\end{lem}
\medskip

\begin{proof}
For every $x \in X, f \in F,$ and $\theta, \varphi \in \mathbb{R}$ we have
$$\sigma(\lambda(f(x) + \theta) + \varphi) 
\left\{ \begin{array}{llll}
\to 1 & \text{for} & f(x)+ \theta > 0 & \text{as} \: \lambda \to +\infty, \\
\to 0 & \text{for} & f(x)+ \theta < 0 & \text{as} \: \lambda \to +\infty, \\
= \sigma(\varphi) & \text{for} & f(x) + \theta = 0 & \text{for all} \:\: \lambda.
\end{array}\right. $$

Then, $\sigma_{\lambda}(x) \overset{\mbox{def.}}{=} \sigma(\lambda(f(x)+\theta)+\varphi)$ converges to
$$\gamma(x) = 
\left\{ \begin{array}{lll}
1 & \text{for} & f(x)+ \theta > 0, \\
0 & \text{for} & f(x)+ \theta < 0, \\
\sigma(\varphi) & \text{for} & f(x) + \theta = 0
\end{array}\right. $$
as $\lambda \to +\infty.$

Now, suppose that $\mu\in B(X)$ satisfies
$$
\int_{X} \sigma(f(x)+\theta) d\mu(x) = 0,\quad\quad\forall f \in F, \, \theta \in \mathbb{R}.$$
Since $\mathbb{R}F\subset F$
(and so $\lambda f$ belongs to $F$, $\forall f\in F$ and $\lambda\in\mathbb{R}$), one has
$$
\int_{X} \sigma_{\lambda}(x) d\mu(x)=0,\quad\quad\forall f\in F,\, \lambda\in\mathbb{R}.
$$
Then,
$$
0 = \lim_{\lambda \to \infty} \int_{X} \sigma_{\lambda}(x) d\mu(x) = \int_{X} \gamma(x) d\mu(x) = \sigma(\varphi) \cdot \mu(\Pi_{f, \theta})+\mu(H_{f, \theta})
$$ 
where $\Pi_{f, \theta} = \{ x : f(x) + \theta =0 \}$ and $H_{f, \theta}= \{ x : f(x) + \theta >0 \}$, for all $f\in F$ and $\theta,\varphi\in\mathbb{R}$.
But $\sigma$ is sigmoidal so by fixing $(f,\theta)\in F\times \mathbb{R}$ and letting $\varphi\to\infty$ above we obtain
$$
\mu(\Pi_{f,\theta})+\mu(H_{f,\theta})=0,\quad\quad\forall f\in F,\, \theta\in \mathbb{R}.
$$
Now, fix $f \in F$ and define the linear functional $G:L^\infty(\mathbb{R})\to\mathbb{R}$ by
$$
G(h) = \int_{X} h(f(x)) d\mu(x).
$$
By taking $\theta\in\mathbb{R}$ and evaluating $G$ at the indicator function $\chi_{[-\theta, \infty)}$ one has
\begin{eqnarray*}
G(\chi_{[-\theta,\infty)})&=&\int_X\chi_{[-\theta,\infty)}(f(x))d\mu(x)\\
&=&\int_{\Pi_{f,\theta}}d\mu+\int_{H_{f,\theta}}d\mu\\
&=&\mu(\Pi_{f,\theta})+\mu(H_{f,\theta})\\
&=& 0,\quad\quad\forall \theta\in\mathbb{R}.
\end{eqnarray*}
Replacing $-\theta$ by $\theta$ in the above result
we obtain
$$
G(\chi_{[\theta,\infty)})=0,\quad\quad\forall \theta\in\mathbb{R}.
$$
It follows that
$$
\int_X\chi_{[\theta,\infty)}(x)d(f_*(\mu))(x)=\int_X\chi_{[\theta,\infty)}(f(x))d\mu(x)=G(\chi_{[\theta,\infty)})=0,\quad\quad\forall \theta\in\mathbb{R}.
$$
By writing $h$ as finite linear combination of indicator functions we obtain that
$$
\int h(x)d(f_*(\mu))(x)=0,\quad\quad\forall h\in L^\infty(\mathbb{R}).
$$
As in \cite{C89}, this proves that Fourier transform of $f_{*} \mu$ equals $0$ for all $f \in F$, and so, $f_{*} \mu = 0 $ for all $f \in F$.
Since $F$ separates $B(X)$, we get $\mu = 0$ and so $\sigma$ is discriminatory with respect to $F$. This completes the proof.
\end{proof}

\begin{proof}[Proof of Proposition \ref{MT2}]
Let $X$ be a compact metric space, and $\sigma : \mathbb{R} \to \mathbb{R}$ be an activation function of sigmoidal type. Suppose that $F \subseteq \mathcal{C}(X)$ separates $B(X)$ and $\mathbb{R}F \subset F$.
Then, $\sigma$ is discriminatory with respect to $F$ by Lemma \ref{lemma2} and so $N_\sigma(F)$ is $C^0$-dense in $\mathcal{C}(X)$ by Lemma \ref{lemma1}.
\end{proof}

Additionally, we include the following lemma regarding the properties of $d_{GH^0}$ as another crucial component.
The proof is omitted due to its similarity to Theorem 1 in \cite{AM17}.

\medskip
\begin{lem}
\label{davis}
The following properties hold for all $f,g,h\in \mathcal{C}$:
\begin{enumerate}
\item
$d_{GH^0}(f,g)=d_{GH^0}(g,f)$.
\item
$d_{GH^0}(f,h)\leq 2(d_{GH^0}(f,g)+d_{GH^0}(g,h))$.
\item
If $D(f)=D(g)$, then
$d_{GH^0}(f,g)\leq d_{C^0}(f,g)$.
\end{enumerate}
\end{lem}
\medskip

Now, we can prove our result.

\begin{proof}[Proof of Theorem \ref{pertinho}]
Let $\mathcal{F}$ be a locally sliceable subset of $\mathcal{C}$ satisfying \eqref{rosa}. Suppose also that
$$
\{X\in D(\mathcal{F}):\mathcal{F}(X)\mbox{ separates }M(X)
\mbox{ and }\mathbb{R}\mathcal{F}(X)\subset \mathcal{F}(X)\}
$$
is dense in $D(\mathcal{F})$.

Let $\sigma:\mathbb{R}\to\mathbb{R}$ be an activation function of sigmoidal type.
Then, Proposition \ref{MT2} applied to $F=\mathcal{F}(X)$ yields that
\begin{equation}
\label{geo}
\{X\in D(\mathcal{F}):N_\sigma(\mathcal{F}(X))\mbox{ is $C^0$-dense in }\mathcal{C}(X)\}
\end{equation}
is dense in $D(\mathcal{F})$.

Now, take $g\in\mathcal{C}$ and $\epsilon>0$.
By \eqref{rosa} there is $\hat{f}\in D^{-1}(D(\mathcal{F}))$ such that
$$
d_{GH^0}(\hat{f},g)\leq\frac{\epsilon}4.
$$
Since $\mathcal{F}$ is locally sliceable, $D(\mathcal{F})$ is a locally sliceable subset of $\mathcal{M}$ by definition.
Take $\Delta=\frac{\epsilon}{16}$ yielding the neighborhood $U$ of $\hat{f}$ as the $\Delta$-ball in $\mathcal{C}$ centered at $\hat{f}$.
For this $U$ we take the $\delta>0$ such that the neighborhood $W$ of $D(\hat{f})$ defined as the $\delta$-ball of $\mathcal{M}$ around $D(\hat{f})$
satisfies the requirements of the definition of locally sliceable fibrewise set. 

Since $D(\hat{f})\in D(\mathcal{F})$, \eqref{geo} provides
$\hat{X}\in D(\mathcal{F})$ such that
\begin{equation}
\label{star}
d_{GH}(\hat{X},D(\hat{f}))<\delta\quad\mbox{ and }\quad
N_\sigma(\mathcal{F}(\hat{X}))\mbox{ is $C^0$-dense in }\mathcal{C}(\hat{X}).
\end{equation}
Then, $\hat{X}\in W$ so
the locally sliceability and the definition of $U$ above provide
$\check{f}\in \mathcal{C}(\hat{X})$ such that
$$
d_{GH^0}(\check{f},\hat{f})<\frac{\epsilon}{16}.
$$
So, the $C^0$-denseness in \eqref{star} provides
$f\in N_\sigma(\mathcal{F}(\hat{X}))$ with
$d_{C^0}(f,\check{f})\leq \frac{\epsilon}{16}$. Since $D(f)=\hat{X}=D(\check{f})$, we can apply Item
(3) of Lemma \ref{davis} to obtain
$$
d_{GH^0}(f,\check{f})<\frac{\epsilon}{16}.
$$
On the one hand,
$f\in N_\sigma(\mathcal{F}(\hat{X}))\subset N_\sigma(\mathcal{F})$
so $f\in N_\sigma(\mathcal{F})$ and, 
on the other,
\begin{eqnarray*}
d_{GH^0}(f,g) & \leq & 2(d_{GH^0}(\hat{f},g)+d_{GH^0}(\hat{f},f))\\
& \leq & 2(d_{GH^0}(\hat{f},g)+2(d_{GH^0}(\check{f},\hat{f})+d_{GH^0}(\check{f},f)))\\
& < & 2 \left(\frac{\epsilon}4+2 \left(\frac{\epsilon}{16}+\frac{\epsilon}{16} \right) \right)\\
& = & \epsilon
\end{eqnarray*}
by items (1) and (2) of Lemma \ref{davis}.
Then, for all 
$g\in \mathcal{C}$ and $\epsilon>0$ we have found $f\in \mathcal{N}_\sigma(\mathcal{F})$ such that
$$
d_{GH^0}(f,g)<\epsilon.
$$
Therefore, $\mathcal{N}_\sigma(\mathcal{F})$ is dense in $\mathcal{C}$ completing the proof.
\end{proof}

\section*{Funding}
WJ was partially supported by the National Research Foundation (NRF) grant by Korea government (MSIT) (No. NRF-2021R1F1A1052631).

\section*{Declaration of competing interest}

\noindent
There is no competing interest.

\section*{Data availability}

\noindent
No data was used for the research described in the article.


\begin{thebibliography}{0}

\bibitem{AM17}
Arbieto, A., Morales C.A.,
Topological stability from Gromov-Hausdorff viewpoint,
{\em Discrete Cont. Dynam. Syst.} 37 (2017) 3531--3544.

\bibitem{B21}
Bueno, C,
Universal Approximation for Neural Nets on Sets,
Thesis (Ph.D.) University of California, Santa Barbara. 2021.


\bibitem{bbi}
Burago, D., Burago, Y., Ivanov, S.,
{\em A Course in Metric Geometry},
Graduate Studies in Mathematics 33, American Mathematical Society, Providence (2001).

\bibitem{5G}
Chazal, F., Cohen-Steiner, D., Guibas, L.J.,  M{\'e}moli, F. Oudot, S.Y., 
Gromov-Hausdorff stable signatures for shapes using persistence, Computer Graphics Forum, 28 (2009),1393--1403.


\bibitem{C89}
Cybenko, G., 
Approximation by superpositions of a sigmoidal function, Mathematics of Control, Signals, and Systems 2 (1989), 303-–314.

\bibitem{ST09}
Ethier, S. N. and Kurtz, T. G.,
{\em Markov processes: characterization and convergence}, John Wiley \& Sons, 2009.



\bibitem{j}
James, I.M.,
{Fibrewise topology},
Cambridge Tracts in Mathematics 91 (1989).



\bibitem{ju}
Jung, W.,
The closure of periodic orbits in the Gromov-Hausdorff space,
{\em Topology Appl.} 264 (2019), 493--497.










\bibitem{lm}
Lee, J., Morales, C.A.,
{\em Gromov-Hausdorff Stability of Dynamical Systems and Applications to PDEs}, Frontiers in Mathematics. Birkh\"auser/Springer,
Cham, (2022). Series ISSN 1660-8046. pp. 160.

















\end{thebibliography}
\end{document}